\documentclass[11pt]{article}

\usepackage{color,graphicx,times}
\usepackage{hyperref}

\newtheorem{definition}{Definition}
\newtheorem{lemma}{Lemma}
\newtheorem{sub-lemma}{Sub-lemma}
\newtheorem{sublemma}{Sub-lemma}
\newtheorem{theorem}{Theorem}
\newtheorem{corollary}{Corollary}
\newtheorem{proposition}{Proposition}

\newtheorem{fact}{Fact}

\newtheorem{observation}{Observation}

\newtheorem{formula}{Formula}
\newcommand{\qed}{\hfill {\bf q.e.d} }
\newcommand{\Bbb}[1]{%
{\bf #1}}
\newcommand{\frak}[1]{%
{\bf #1}}
\newcommand{\proof}{ {\bf Proof:}}
\newcommand{\sech}{ {\tt sech}}

\input epsf.sty

\begin{document}

\title{The ideal Thurston-Andreev theorem and triangulation production}
\author{Gregory Leibon \\ Department of Mathematics  \\  
Dartmouth College\\ Hanover, NH}
\maketitle


\vspace{1in}
\begin{center}
{ mailing address:}
\end{center}


\begin{center}
{\large{Gregory Leibon}} \\
{\large 29 Fletcher Circle} \\
{\large Hanover, NH} \\
{\large 03755}
\end{center}

\vspace{1in}

\begin{center}
{ e-mail address:}
\end{center}


\begin{center}
{\large gleibon@dartmouth.edu}
\end{center}
\pagebreak

\begin{section}{Introduction}
The main result in this paper is a generalization 
of the convex ideal case 
of the Thurston-Andreev theorem when $\chi(M) <0$. 
The proof naturally decomposes into its non-linear and linear aspects.
The non-linear part  is essentially  a 
triangulation  production theorem,  stated in section \ref{san}.  
This theorem concerns taking a 
topological triangulation  along with  
formal angle data and  
then ``conformally 
flowing'' this formal angle data
to uniquely associated uniform angle data, where uniform angle data means  
the data contained in a  geodesic 
triangulation of a hyperbolic surface.
This flow is the gradient of an objective function related in a rather
magical way to
hyperbolic volume.   That such a magical connection might exist was first
explored in Bragger \cite{Be}, and the hyperbolic volume needed in the
case presented here was observed by my thesis advisor, Peter Doyle.

Conformally flowing turns out to be related to to certain disk patterns and
hyperbolic polyhedra, as discussed in section \ref{poly} and \cite{Le}. The
linear part of this paper  concerns  gaining an explicit handle  on 
which patters can arise.
In the end, a complete characterization of
the ``convex ideal" patterns in the  $\chi(M)<0$ 
case of the Thurston-Andreev theorem is presented,  for the statement of
which see section \ref{sttur}.  It worth noting that in the torus  and 
spherical cases that this generalization had already been 
accomplished. 
The toroidal case of the
entire strategy used here    
has its origins 
in the beautiful and often overlooked work of 
Bragger \cite{Be} 
and can also be found  
in  \cite{Ri2}. The  spherical  case 
of this generalization of the Andreev theorem
was accomplished by Rivin in \cite{Ri1}.
Whether there is proof 
directly using  the spherical version of the techniques 
in this paper  is still unknown.    
For an account of the theorem being  generalized
see  Thurston's \cite{Th}.

This paper is organized as follows: section
\ref{tdt} contains the statements, set up, and notation needed to describe the
theorems mentioned above. 
Section \ref{prot} contains  the 
proof of the triangulation production theorem using a bit
of hyperbolic geometry.  
Section \ref{linn} contains the proof of its corollary, the
generalization of the Thurston-Andreev theorem stated in section
\ref{sttur}; in the form of a  
``min flow max cut'' type argument.  In the final section a discussion of
some known generalizations and some questions takes place.


I would like to thank my thesis advisor 
Peter Doyle for sharing his many beautiful ideas with me; 
without him the work here would 
not have been possible.

\end{section}

\begin{section}{Statements and Notation}\label{tdt}
\begin{subsection}{The Triangular Decomposition Theorem }\label{san}
Throughout this paper $M$ will  denote a compact two-dimensional 
surface 
with $\chi(M) < 0$.  
By geometry I   will mean a hyperbolic structure.  
Uniqueness of geometries, triangulations and disk patterns is 
of course up to isometry.  

The main theorem in this section  really 
should be stated for the following 
structure, which generalizes the notion of triangulation.

\begin{definition}\label{decomp}
Let a  triangular decomposition, $\frak{T}$, 
be a cell decomposition of $M$  that  lifts to a triangulation
in $M$'s universal cover. 
\end{definition} 

We will 
keep track of the combinatorics of such a decomposition by denoting 
the vertices as $\{ v_i \}_{i=1}^{V}$, 
the edges  as $\{ e_i \}_{i=1}^{E}$  and the triangles 
as $\{ t_i \}_{i=1}^{F} $. Let $E$, $V$, $F$, $\partial
E$ and $\partial V$  denote the sets of edges, vertices, faces, boundary
edges, and boundary vertices respectively.  For convenience this same
notation will denote the cardinalities  of these sets. Let
$\{e \in S\}$  denote   the set of   edges on the surface in a collection  of 
triangles $S$,
and let $\{e \in v\}$  denote all the edges associated to a vertex $v$ 
as if counted in the universal cover.  The set of triangles containing a
vertex  $\{ t \in v\}$ has the special name of the flower at $v$.

\begin{figure*}
\vspace{.01in}
\hspace*{\fill}
\epsfysize = 1.5in  
\epsfbox{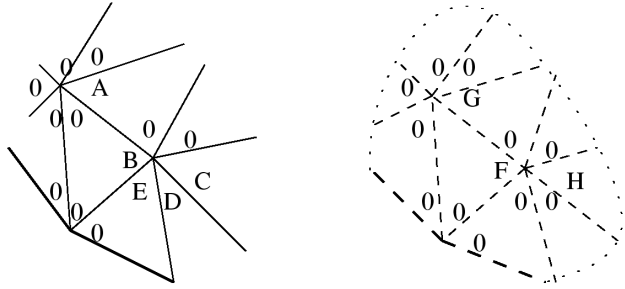}
\hspace*{\fill}
\vspace{.01in}
\caption{\label{the5} {\bf A covector and a vector in geometric notation.} 
A vector will be denoted
by placing its  coefficients, $\{A^i\}$, in a copy of the triangulation 
with dashed lines and covector will contain its  
coefficients, $\{A_i\}$, 
in a copy of the triangular decomposition with solid lines.  
Thick lines will
always denote a boundary edge, as  in lower the left corner of 
pictured vector and covector.
If in the picture we mean the non-specified values 
to be arbitrary we will surround the picture with a loop (as with the pictured
vector)  and  if we mean   the non-specified values  to be zero 
the picture will not be surrounded (as with the pictured covector). }
\end{figure*}

Note that in a triangular decomposition there are 
$3F$ slots $\{ \alpha_i \}$ in which one can insert 
possible triangle angles, which we 
will place an order on and identify with a basis of a 
$3 F$ dimensional real vector space.  With this basis choice we will 
denote this vector space as   
$\Bbb R^{3F}$, and  denote vectors in it as $x =\sum A^i \alpha_i$.
Further more let 
$\alpha^i$ be a dual vector such that $\alpha^i(\alpha_j) = \delta_i^j$. 
With this we will view the 
angle at the slot $\alpha_i$ as  $\alpha^i(x) = A^i$.
It is rarely necessary to use this notation and instead 
to use the actual geometry as in figure \ref{the5}. 
Notice the pairing of a vector and a covector 
denoted 
can be viewed geometrically as in figure \ref{the6}. 
For a  triangle $t$ containing the angle slots 
$\alpha_i$, $\alpha_j$, and $\alpha_k$ 
let $d^t (x)   = \{ A^i,A^j,A^k \}$
and call $d^t(x)$  the angle data associated to $t$. 

\begin{figure*}
\vspace{.01in}
\hspace*{\fill}
\epsfysize = 1.5in  
\epsfbox{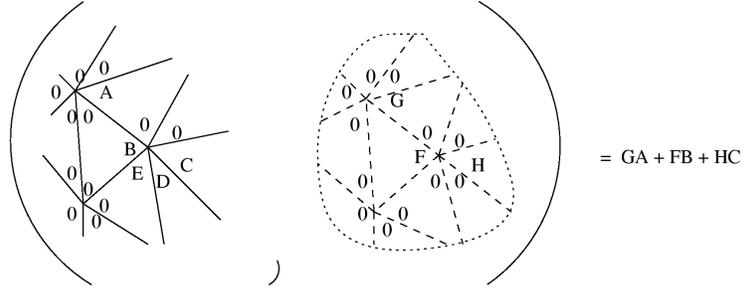}
\hspace*{\fill}
\vspace{.01in}
\caption{\label{the6} { {\bf The geometric pairing of a vector and a covector}
is achieved  by placing the copy of the   triangular decomposition corresponding
to the vector on top of the  triangular decomposition  corresponding to the
covector and multiplying  the numbers living in the same angle slots.}}
\end{figure*}

\begin{figure*}
\vspace{.01in}
\hspace*{\fill}
\epsfysize = 1.5in  
\epsfbox{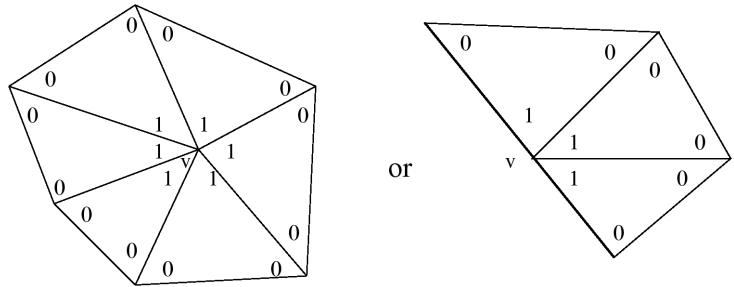}
\hspace*{\fill}
\vspace{.01in}
\caption{\label{co0} For each vertex $v$ the covector indicated in this
figure will be denoted $p^v$}
\end{figure*}

In order to live on an actual nonsingular geometric surface
with geodesic boundary all such angles should be required to live
  in the subset of 
$\Bbb R^{3F}$ where the angles at an interior 
 vertex sum to $2 \pi$ and the angles at a boundary vertex sum to 
$\pi$. 
Using the in figure \ref{co0},
this encourages us to choose  our possible angles in the affine flat
\[ \frak{V} = \{x \in \Bbb R^{3F} \mid p^v(x) = 2 \pi \mbox{ for all } v 
\in V -\partial V \mbox{ and } p^v(x) =  \pi 
\mbox{ for all } v \in \partial V  \}. \]

To further limit down the possible angle values we define 
the covector $l^t$ as in figure \ref{co2}
\begin{figure*}
\vspace{.01in}
\hspace*{\fill}
\epsfysize = 1in  
\epsfbox{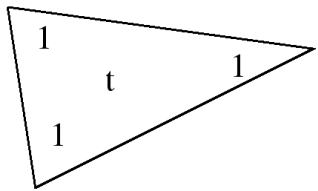}
\hspace*{\fill}
\vspace{.01in}
\caption{\label{co2} For each $t$ we will denote the pictured covector as
$\sigma^t$.}
\end{figure*}
Using the notation in figure \ref{co2}, by the Gauss-Bonnet theorem we know that
$k^t(x)$ defined as 
\[ k^t(x) = \sigma^t(x) - \pi \] 
would be the curvature in  a geodesic triangle with 
angle data $d^t(x)$. 
We will now isolate the open convex subset of $\frak{V}$ 
where the curvature is negative and all angles are realistic.
\begin{definition}\label{dcc}
Let an angle system be a point in  
  \[ \frak{N} = \{ x \in \frak{V} \mid k^t(x) < 0 \mbox{ for all } t 
 \mbox{ and  } \alpha^i(x) \in (0, \pi) \mbox{ for all } \alpha^i\}. \]
\end{definition}
Note the actual angle data of a geodesic 
triangulation of a  
surface with negative curvature  
has its angle data living in this set.

Observe  if $k^t(x) < 0$ then we may form an actual hyperbolic triangle with
the 
angles in $d^t(x)$.   For each $e \in t$ denote
the length of the edge $e$ with respect to this data as
$l^e_t(x)$.

Suppose the triangles of $\frak{T}$  fit together in the sense that 
$l^e_{t_1}(x) = l^e_{t_2}(x)$ whenever it makes sense.  Then
$\frak{T}$ being a triangular decomposition implies every open 
flower 
is embedded in $M$'s universal cover and  
when the edge lengths all agree
this flower can  be given  
a hyperbolic structure which is consistent on flower overlaps.  
So we have formed  a hyperbolic structure
on $M$. 

\begin{definition}\label{duas}
Call  an angle system $u$ uniform 
if all the hyperbolic realizations of the triangles in $u$
fit together to form a hyperbolic structure on $M$. 
\end{definition}

In section 2.4  we will  attempt to take a point in $\frak{N}$ and 
deform it into a uniform point.  Such deformations are located in  an
affine space and I will call them conformal deformations 
(see \cite{Le}  for a more careful motivation of this terminology).  
To describe this  affine space for each edge $e \in E -\partial E$ 
construct a vector $w_e$ as in figure \ref{ve1}.
\begin{figure*}
\vspace{.01in}
\hspace*{\fill}
\epsfysize = 1in  
\epsfbox{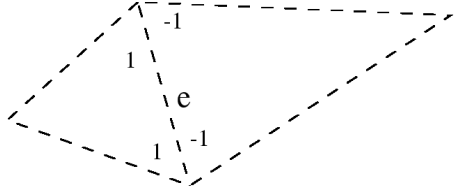}
\hspace*{\fill}
\vspace{.01in}
\caption{\label{ve1} For each edge $e \in E -\partial E$ 
let $w_e$ denote a vector as in this figure.}
\end{figure*} 
\begin{definition}\label{condef}  
Using the notation from figure \ref{ve1}, a conformal deformation will be a
vector in 
\[ C = span\{ w_e \mid \mbox{ for all } e \in E - \partial E \}.\]
Call
$x$ and $y$ conformally equivalent
if $x - y \in C$, and
let 
\[\frak{N}_x = (x + C) \bigcap \frak{N}\]
denote the  conformal class of $x$. 
\end{definition}
The first thing worth noting is that if $x \in \frak{V}$ and 
$y$ is conformally equivalent to $x$, then looking at the  
pairing between  $p^v$ with $v \in V - \partial V$ and the 
$w_e$ (as in
figure
\ref{co0}) we have  
\[ p^v(y) = p^v\left(x + \sum_{e \in \frak{P}} B^{e} w_{e}\right) 
= p^v(x) + \sum_{e \in \frak{P}} B^{e} p^v(w_{e} ) 
= 2 \pi  + 0, \]
hence $y$ is also in $\frak{V}$.    Similarly for $v \in \partial V$.

To combinatorially understand the points in $\frak{N}$ 
which we may  conformally deform into uniform 
structures it is useful to express a particularly 
nasty set in the boundary of $\frak{N}$.  
\begin{definition}\label{legal}
Let $t$  be called a legal with respect to $x \in \partial \frak{N}$ if 
$d^t(x)  = \{A^1,A^2,A^3\}  \neq \{ 0,0,\pi \}$ yet either $k^t(x) =0$
or for some $i$ we 
have $A^i =0$. Let  
\[ \frak{B} = \{x \in \partial \frak{N} \mid  \mbox{ x contains no legal  } t
\}.
\] 
\end{definition}

We will prove the following theorem.
\begin{theorem}\label{tri1}
If there is a uniform angle system conformally 
equivalent to $x$ then it is unique, and  for any angle 
system $x$ with   $(x + C) \bigcap B$ empty 
there exists a conformally equivalent uniform angle system. 
\end{theorem}

Much of what takes place here relies on certain basic 
invariants of conformal 
deformations, which is the subject of the next section.

 \end{subsection}

\begin{subsection}{Ideal Disk Patterns}\label{poly}
In this section we introduce the disk patterns that a
conformal class is related to. 
We begin with the
introduction and interpretation of certain key conformal
invariants.    Using the notation in figure \ref{co3}
for each edge $e \in \partial V$  we will denote  $\psi^e_t$ as $\psi^e$
while  for each edge $e \in E - \partial E$ associated with triangles
$t_1$ and $t_2$ we will let
\[ \psi^{e} =  \psi^{e_1}_{t_1} + \psi^{e_1}_{t_2}.\]    
We will call the  $\psi^e$ covector  the formal angle complement at $e$.
Let the formal intersection angle be defined as
\[ \theta^e(x) = \pi - \psi^e(x)\]
when $e \in E -\partial E$.
 and 
\[ \theta^e(x) = \frac{\pi}{2} - \psi^e(x)\]
when $e \in \partial E$.

\begin{figure*}
\vspace{.01in}
\hspace*{\fill}
\epsfysize = 1in  
\epsfbox{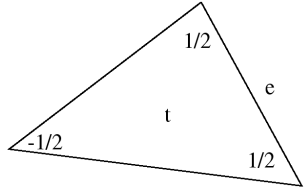}
\hspace*{\fill}
\vspace{.01in}
\caption{\label{co3} For each triangle $t$ and $e \in t$ we will let 
$\psi^e_t$ denote the pictured covector.}
\end{figure*} 
Looking at the pairing 
between a  $\psi^e$ and  $w_e$ (as in figure \ref{the6}),
we see if $y$ is conformally equivalent to $x$ then 
\[ \psi^e(y) = \psi^e\left(x + \sum_{f \in \frak{T}} B^f w_{f}\right) 
= \psi^e(x) + \sum_{f \in \frak{T}} B^f \psi^e(w_{f} ) 
= \psi^e(x),\]
and indeed for each relevant edge $e$ we see 
$\psi^e$ and $\theta^e$ are conformal invariants.

The  fundamental 
reason why theorem \ref{tri1} is  
related to circle patterns and polyhedra construction is the following
wonderful observation.

\begin{observation}\label{fund}
At a uniform angle system $u$  and 
$e \in E -\partial E $,  $\theta^e(u)$ is 
the intersection angle
of the circumscribing circles of the two 
hyperbolic triangles sharing $e$.
\end{observation}

\proof
The proof given here relies on a principle fundamental to everything 
that takes place here,  namely that questions taking place on a hyperbolic
surface can often be best interpreted by viewing the question in
three-dimensional  hyperbolic space,  $H^3$.   Let us call $H^2$ the
the hyperbolic plane in $H^3$ viewed as  in figure
\ref{hit}.  
The inversion, $I$,  through the sphere of radius $\sqrt{2}$
centered at the south pole interchanges our specified $H^2$ with the upper
half of  the sphere at infinity $S^{\infty}_u$.  Notice when viewed
geometrically  this map sends a point
$p \in H^2$ to the point where the geodesic perpendicular to $H^2$
containing $p$   hits
$S^{\infty}_u$ (see figure \ref{hit}).  In particular being an inversion 
any circle in the $xy$-plane is mapped to a circle on the sphere at $\infty$,
$S^{\infty}$.
\begin{figure*}
\vspace{.01in}
\hspace*{\fill}
\epsfysize = 2in 
\epsfbox{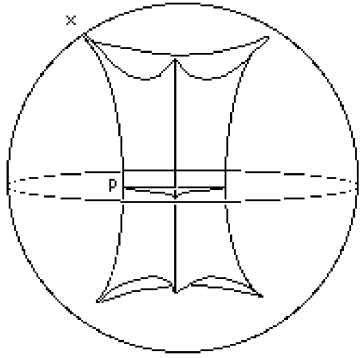}
\hspace*{\fill}
\vspace{.01in}
\caption{\label{hit} In the figure  we have our specified hyperbolic plane
$ H^2 \subset H^3$  realized as the intersection of the unit sphere at the
origin with the $xy$-plane in $\Bbb{R}^3$ via the Poincare  disk model of
$H^3$.  The point $p$ in the figure  is mapped  to the point
labeled $x$ under the inversion $I$.   We also are viewing a triangle on
that plane and its associated prism in this model.}
\end{figure*}

The use of this mapping will require the  introduction of an object that will
be  crucial in proving theorem  \ref{tri1}.

\begin{definition} \label{pridef}
Place a hyperbolic triangle on a copy of $H^2 \subset H^3$. 
Let its associated prism be the convex hull of 
the  set  consisting of the triangle 
unioned with  the geodesics 
perpendicular to this $H^2 \subset H^3$ 
going through the triangles
vertices as visualized in figure \ref{hit}.  
Denote the prism  relative to the hyperbolic triangle
constructed from the data $d = \{A,B,C\}$ as 
$P(d)$.
\end{definition}

 Now back to our proof.  Let $u$ be a uniform angle system and let $t_1$
and $t_2$ be a pair of triangles sharing the edge $e$.  Place them next to
each other in our $H^2$ from figure \ref{hit}.  Notice
$t_1$ and $t_2$ have circumscribing  circles in the $xy$-plane, which
correspond to either circles, horocirles, or bananas in the $H^2$ geometry. 
Since the Poincare model is conformal the intersection angle of these
circles is precisely the hyperbolic intersection angle.   
Being an inversion $I$ is conformal, so these circles are sent
to    circles at infinity intersecting at the same angle  and going
through the ideal points of the neighboring $P(d_{t_1}(u))$ and
$P(d_{t_2}(u))$.   But these  circles at infinity  are also the
intersection  of $S^{\infty}$ with the spheres representing 
the hyperbolic planes forming the top faces of 
$P(d_{t_1}(u))$ and $P(d_{t_2}(u))$.   
So the intersection angle of these spheres is precisely the sum
of the angles inside $P(d_{t_1}(u))$ and $P(d_{t_2}(u))$ at the edge
corresponding to $e$,   which we will now see is $\theta^e(u)$.  
In fact we will show that this decomposition of the intersection angle
is precisely the decomposition
 \[ \theta^e(u)= \theta^e_{t_1}(u) + \theta^e_{t_2}(u).\]

\begin{figure*}
\vspace{.01in}
\hspace*{\fill}
\epsfysize = 2.5in 
\epsfbox{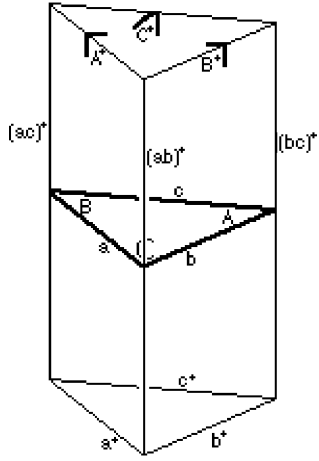}
\hspace*{\fill}
\vspace{.01in}
\caption{\label{clam} The notation for an ideal prism associated to 
hyperbolic angle data $\{A,B,C\}$, viewed for convenience  in the Klein
model.  }
\end{figure*}

Now pick an $i$ and let  $d^{t_i}(u)$ be denoted  by $\{A,B,C\}$.  Assume
our specified edge
$e$ corresponds to the
$a$ in figure \ref{clam}.
From figure \ref{horo1} we see the angles in figure \ref{clam} satisfy the
system of linear equation telling us that interior angles of the prism sum
to
$\pi$ at each vertex of the prism.  Solving this system   for the needed
angle, 
$A^{\star}$, we find that indeed
\[ A^{\star} =\frac{\pi + A -B -C}{2 } = \theta^e_{t_i}(u), \]
as claimed.

\qed

\begin{figure*}
\vspace{.01in}
\hspace*{\fill}
\epsfysize = 2in 
\epsfbox{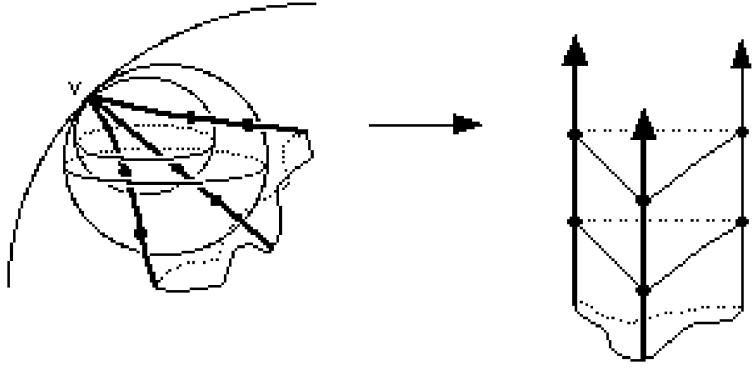}
\hspace*{\fill}
\vspace{.01in}
\caption{\label{horo1} In this figure we have  cut off an ideal
vertex of a convex hyperbolic polyhedron with a horoball. Horospheres are
flat planes in $H^3$ and our cut produces an Euclidean polygon in
this flat plane.  In particular the interior angles at an ideal vertex of a
convex hyperbolic  polyhedron are those of an Euclidean polygon.  As
indicated, this is best  seen the upper half space model of
$H^3$, with the ideal vertex being sent to
$\infty$.  (Nice applications of this observation can be found in \cite{Th}.) }
\end{figure*}

At this point it is useful to name the appropriate 
 home of the possible intersection angle assignments.
Just as we did with the angle 
slots, let the edges correspond to the basis vectors of 
an $E$  dimensional vector space, which we will denote 
$\Bbb R^E$ with this 
basis choice.  We will be viewing this as the space of 
possible angle complements. Denote these vectors as 
$p = \sum_{e \in \frak{T}} E^{e} \psi_e$.
let
\[ \Psi: \Bbb R^{3 F} \rightarrow \Bbb R^{E} \]
be the linear mapping 
given by 
\begin{eqnarray}\label{mapping}
 \Psi(x) = \sum_{e \in \frak{T}} \psi^{e}(x) \psi_e,
\end{eqnarray}
and note  
by observation \ref{fund} that we do indeed hit the intersection  angle
discrepancies when applying $\Psi$ to an uniform angle system.
Let $\psi^e$ denote the algebraic dual to $\psi_e$.  To justify this
abuse of notation, note that 
$\psi^e(\Psi(x)) = \psi^e(x)$. Similarly let
$\theta^{e}(p) = \pi - \psi^e(p)$ for
$e \in E -\partial E$ and 
$\theta^{e}(p) = \frac{\pi}{2} - \psi^e(p)$
 when 
$e \in \partial E$.

\begin{definition}
A point $p \in \Bbb{R}^{E}$ relative to a triangular
decomposition $\frak{T}$ will be called   an ideal disk pattern if $p =
\Psi(u)$ for some uniform angle system $u$.

\end{definition}

In less dramatic terms this simply says that such a $p$ can be realized
geometrically.   Theorem \ref{tri1}
provides us with the following corollary telling us such realizations are
unique.

\begin{corollary}\label{idpun}
For any $p\in \Bbb{R}^E$ at most one uniform $u$ can satisfy  $\Psi(u)
= p$.   
\end{corollary}
{\bf Proof:}
$\Psi$ has rank $E$ since the pairing of $\psi^{e}$ with 
the vector $m_{e}$  in figure \ref{ve2} 
satisfies   $\Psi(m_{e}) = \psi_e$ for each edge $e$.  By the conformal
invariance noted in the previous section  the null space contains the 
$E  - \partial E$ dimension space $C$ and is 
\[ 3F - E = 2E - \partial E - E= E -\partial E\]
dimensional,  so $C$ is precisely the  null space.
In particular all angle systems 
which could conceivably hit a specified $p$ will be in  
$\Psi^{-1}(p)$, which is $x +C $ for some $x$.
So theorem \ref{tri1} guarantees the uniqueness of the associated angle
system. 

\qed

\begin{figure*}
\vspace{.01in}
\hspace*{\fill}
\epsfysize = 1in  
\epsfbox{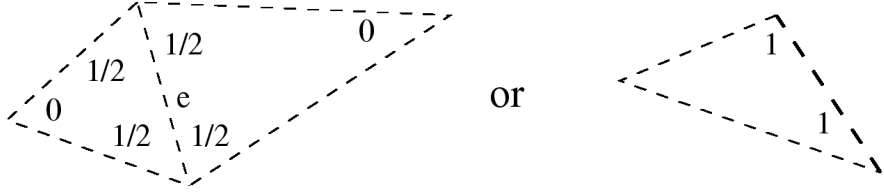}
\hspace*{\fill}
\vspace{.01in}
\caption{\label{ve2} For each edge let  $m_e$ be the pictured vector.}
\end{figure*}

In the next section we will derive affine conditions 
on a point $p \in \Bbb{R}^E$ relative to a triangular decomposition
$\frak{T}$ to be an ideal disk pattern.

I will finish this section by noting that such patterns are related to
a certain family of  hyperbolic polyhedra.
Namely if we take the polyhedra constructed by placing
a geodesic triangulation on an $H^2 \subset H^3$ and forming  
$\bigcup_{t \in \tilde{\frak{T}}} P(d_t(u))$.  
Note as a schollium to fact \ref{fund} that the dihedral angles in the
polyhedra are precisely the intersection angles in the disk patterns.   So
any question about such ideal disk patterns can be translated into a
question concerning such polyhedra (see \cite{Le} for details).   In
particular in the following sections we will be placing constraints on the
$\theta^e$ to be in $(0,\pi]$,  and its worth noting
that this corresponds precisely to the convex case among these
polyhedra.   The results there can be viewed as a characterization
of the possible dihedral angles in such convex polyhedra.

\end{subsection}
\begin{subsection}{A Thurston-Andreev Type Theorem}\label{sttur}
In the previous section we learned that an ideal disk pattern
is always unique when it exists, and we are now left to deal with the 
dilemma of finding good existence criteria.  
This  will be accomplished by noting the restrictions placed on
$ \Psi(x)$ for $x \in \frak{N}$.  For example  the
curvature assumptions in
$\frak{N}$  produces the following restrictions:

\begin{sublemma}\label{ob}
When $x \in \frak{N}$ we have
$\psi_{t}^e(x) \in \left(\frac{-\pi}{2},\frac{\pi}{2} \right)$.
\end{sublemma}
{ \bf Proof:}
Let  $d^{t}(x) = \{A,B,C\} $ and note since 
$B +C \leq A+B+C = \sigma^t(x) < \pi$ and $A<\pi$ we have  
\[ -\frac{\pi}{2} < -\frac{A}{2} <  
\psi^e_{t} = \frac{B+C-A}{2}
<  \frac{B +C}{2} < \frac{\pi}{2} ,\]
as needed.

\qed

In this section we will be strengthening this restriction to the strict convex
case when the ideal disk pattern satisfies
$p \in (0,\pi)^{(E - \partial E)}  \times (0,\frac{\pi}{2})^{\partial
E}$.  An angle system producing such an ideal disk pattern will live in
\[\frak{D} = \left\{ x \in \frak{N} \mid  
\Psi(x) \in  (0,\pi)^{(E - \partial E)}  
\times \left(0,\frac{\pi}{2} \right)^{\partial E} \right\}, \]
and let us call set this the set of
Delaunay angle systems. 
Aside from arising naturally in practice from the 
Delaunay triangulation
of a set of points (see \cite{Le}), these angle systems  
are remarkably  easy to work with because they completely
eliminate the possibility of a conformally equivalent badly
behaved sets of angles.

\begin{lemma}\label{thing}
$x \in \frak{D}$ is never conformally equivalent to 
a point in $\partial \frak{N}$
where for some triangle
$d^t(x) = \{0,0,\pi\}$.
\end{lemma}
\proof
To see this  fact assume to 
the contrary that for some $t$  and $c$ we have 
$d^t(x + c) = \{0,0,\pi\}$. 
Let $e$ be the  edge of $t$ across from
$t$'s  $\pi$
and let $t_1$ be $t$'s neighbor next to $e$ if it exists.

Note by sublemma \ref{ob}
that the conformally invariant $\psi^e(x) \in (0,\pi)$
 would (even 
in the best possible case when $e$ is not on the boundary) 
have 
to satisfy
the inequality 
$\psi^e(x+C) = - \frac{\pi}{2} + \psi_{t_1}^e \leq 0 $,
contradicting the fact $x \in \frak{D}$.

\qed

The elimination of such possibilities allows us to immediately apply theorem 
\ref{tri1}, and arrive at...

\begin{corollary}\label{delnice}
Every point of $\frak{D}$ has a unique ideal 
disk pattern associated to it.
\end{corollary}
{\bf Proof:} 
Let $x \in \frak{D}$ and note we are searching for a uniform
angle system in 
 $ x+C \bigcap \frak{N} = \Psi^{-1}(\Psi(x)) \bigcap \frak{N}$.  From lemma
\ref{thing} we have  that if $x \in \frak{D}$ then no element in 
$x +C$ could possibly be in $B$ and the 
corollary  follows from
 theorem   
\ref{tri1}.

\qed

The goal at this point is to determine necessary and sufficient 
conditions on a point in $\Bbb{R}^{E}$
 to insure that it is $\Psi(x)$ for some Delaunay angle system.  
We immediately have that any such point is $p \in (0,\pi)^{(E -
\partial E)} 
\times (0,\frac{\pi}{2})^{\partial E}$, and our first  non-trivial 
necessary condition  is the condition 
related to the fact that the angles at the internal 
vertex in a geometric 
triangulation sum to $2 \pi$ and at a  boundary vertex sum to $\pi$.

\[ \mbox{(}n_1\mbox{)        } \left\{ \begin{array}{ll}
\sum_{e \in v} \psi^e(p) = 
2 \pi   &   \mbox{  if   } v  \in V - \partial V \\
 \sum_{e \in v} \psi^e(p) = 
\pi   &   \mbox{ if  }  v \in \partial V
\end{array} \right. \]

\proof
To see the necessity of $(n_1)$  we will show
\[ \Psi(V) = \{p \in \Bbb R^{E} \mid p \mbox{ satisfies } (n_1)\}.\]

First note that if $p = \Psi(x)$  then 
\[ \sum_{e \in v} \psi^e(p) =  \sum_{e \in v} \psi^{e}(x) =
p^{v} (x).  \]
So by choosing $x \in V$ we see $\Psi(V)$ is included in 
\[ W = \{p \in \Bbb R^{E} \mid p \mbox{ satisfies } (n_1)\} .\]
Recall from the proof of corollary \ref{idpun} that
$\Psi(\Bbb R^{3F}) = \Bbb R^E$. So 
we may express any $p \in W$ as $p = \Psi(x)$ 
and  the above computation  guarantees $x \in V$ as needed. 

\qed

The second necessary  condition is a global one; namely 
an  insistence that  for every set $S$ of
$|S|$  triangles  in $\frak{T}$ that 
\begin{eqnarray*}
(n_2)  & &   \sum_{e \in S}  \theta^e(p) > \pi |S|.  \\
\end{eqnarray*}

\proof
Verifying  $(n_2)$'s necessity relies on the following formula.
\begin{formula}\label{sett}
Given a set of triangles $S$ 
\[ \sum_{\{e \in S \}} \theta^e(x)   = 
\sum_{t \in S} \left( \pi  - \frac{k^t(x)}{2}\right) 
+ \sum_{e \in \partial S - \partial E}  
\left(\frac{\pi}{2} - \psi^e_{t}(x)\right) , \]
with the   $t$ in $\psi^e_{t}(x)$ term being the triangle on 
the  non-$S$ side of $e$.   
\end{formula}

\[ \sum_{\{e \in S \}} \theta^e(x) = 
\sum_{e \in S - \partial E} (\pi - \psi^e(x)) 
+ \sum_{e \in S \bigcap \partial E} 
\left(\frac{\pi}{2} - \psi^e(x) \right) \]
\[ =
\sum_{e \in S - \partial E} 
\left( \left(\frac{\pi}{2} - \psi_{t_1}^e(x)\right) 
+ \left(\frac{\pi}{2} - \psi_{t_2}^e(x)\right)\right) 
+ \sum_{e \in S \bigcap \partial E} 
\left(\frac{\pi}{2} - \psi^e(x) \right)\] 
\[ = \sum_{t \in S} \left( \pi  + \frac{ \pi - \sigma^t(x)}{2}\right) 
+ \sum_{e \in \partial S - \partial E}  
\left(\frac{\pi}{2} - \psi^e_{t}(x)\right)  \]
with the   $t$ in $\psi^e_{t}(x)$ term being the triangle on 
the  non-$S$ side of $e$.   Substituting the definition of 
$k^t(x)$ gives the needed formula.

To apply this formula note for any point $x \in \frak{N} $ that 
 $- k^t(x) > 0$ and from  sublemma \ref{ob} 
 that $\frac{\pi}{2} - \psi^e_{t}(x) >0$.  
So removing these terms from the above formula strictly 
reduces its size and when summed up we arrive at $(n_2)$.

\qed 

With these two necessary conditions 
 we arrive at a pattern existence theorem 
(see section \ref{itat} for a
stronger result):

\begin{theorem}\label{cir11}
If 
\[ p \in D  = \left\{q \in (0,\pi)^{(E - \partial E)}  
\times \left(0,\frac{\pi}{2}\right)^{\partial E} \mid  q \mbox{ satisfies } 
(n_1) \mbox{ and } (n_2) \right\} \] 
then $p$ is  an 
ideal disk pattern.
\end{theorem}

By corollary  \ref{delnice} above 
and the fact that $(n_1)$  and $(n_2)$ are necessary  
this theorem would follow if 
we knew 
\[\Psi(\frak{D}) = D.\]

Notice the fact that $(n_1)$ and $(n_2)$ are 
necessary guarantees that $\Psi(\frak{D}) \subset D$, 
and we left 
to explore  $\Psi$'s surjectivity.  
It is this bit of linear algebra and will be dealt with in section
\ref{linn}.

\end{subsection}

\end{section}

\begin{section}{The Non-linear Argument}\label{prot}
The non-linear argument is dependent on a rather remarkable link between
hyperbolic volume and uniform structures.
Let the volume of the prism $P(d^t(y))$ be denoted
$V(d^t(y))$, we will be exploring the objective function  
\[H(y) = \sum_{t \in \frak{T}} V(d^t(y)) \]
on $\frak{N}_x$.  First let us note this function has the correct objective.

\begin{fact} 
$y$ is a critical point of $H$ on $\frak{N}_x$ if and only if $y$ is uniform.
\end{fact}
\proof  
To see this fact requires an understanding of 
$H$'s differential at each $x$. 
As usual for 
a function on a linear space like $\Bbb R^{3F}$
we  use translation to 
identify the tangent and cotangent spaces at
every point with  $\Bbb R^{3F}$ 
and  $(\Bbb R^{3F})^*$, and express our differentials in the chosen basis.

 \begin{lemma} \label{diff}
$dH(z) = \sum_{t \in \frak{T}} \left( \sum_{e \in t} h^e_t(z) d \theta^e_{t}
\right)$ with the property that
$h^e_{t}(z)$ uniquely determines the length $l^e_{t}(z)$.
\end{lemma}

This formula will be proved in section \ref{diff}, for now lets see how to use
it.   Recall that 
$T_z(\frak{N}_x) = span \{w_e \mid e \in E -\partial E\}$  from definition
\ref{condef}.
Now simply observe that $\theta^e_t(w_f) = \pm \delta^e_f$ with the sign
depending on whether $t$ contains the negative of positive half of $w_e$. So at
a critical point $y$ we have  
\[ 0  =  dH(y)(w_e) = h^e_{t_1}(y) - h^e_{t_2}(y) \]
where $t_1$ and $t_2$ are the faces sharing $e$.  So from the above lemma
and the fact that the $w_e$  span $T_y(\frak{N}_x)$
 we have that
$l^e_{t_1}(y) = l^e_{t_2}(y)$ for each edge is equivalent
to  $y$ being critical, as needed. 

\qed

We may now prove the  uniqueness
assertion  in  theorem \ref{tri1}.

\begin{proposition}\label{lemer}
If $\frak{N}_x$ contains a uniform angle system this angle system is unique in
$\frak{N}_x$. 
\end{proposition}
{\bf Proof:}
The following lemma will be proved in section \ref{con}. 

\begin{lemma}\label{cony}
$H$   is a strictly concave,  smooth function 
on $\frak{N}_x$ and continuous on $\overline{\frak{N}}_x$.
\end{lemma}

Now a smooth strictly concave function like $E$  has a most one critical point in any
open convex set, which  proves proposition \ref{lemer}.

\qed

Now its time to explore the existence  of critical points.
Given a pre-compact open set $U$ 
and a 
continuous function $F$ on $\overline{U}$
we automatically achieve a maximum. 
For this maximum to be a  critical point it is enough to know that 
$F$ is differentiable in $U$ and 
that the point of maximal $F$ is in
$U$.   

One way to achieve this is to show that for any boundary point $y_0$ 
that there is a direction $v$,
 an $\epsilon >0 $ and a $c >0$  such that $l(s) = y_0 + s v$ satisfies 
\[  l(0,\epsilon) \subset U\] 
and  
\[ \lim_{s \rightarrow 0^+} \frac{d}{ds}F(l(s)) \mbox{  } > \mbox{  } c, \]
for all $s \in (0,\epsilon)$.  
This follows since under these hypotheses   $F(l(s))$ is continuous and 
increasing on $[0,\epsilon)$ and $y_0$ certainly could not have 
been a point where $F$ achieved its maximum.

In our setting we have the following lemma to be proved in section
\ref{bound}.

\begin{lemma}\label{pushin}
For every  point $y_0$ in $\partial \frak{N}$  
and every direction $v$ 
such that 
\[  l(0,\infty) \bigcap  \frak{N} \neq \phi\] 
and
\[  l[0,\infty) \bigcap  \frak{B} = \phi   \]  
we have  
\[ \lim_{s \rightarrow 0^+} \frac{d}{ds}H(l(s))   
\mbox{  } =   \mbox{  }   \infty. \]
\end{lemma}

By convexity of $\frak{N_x}$ for each boundary point there is such a $v$, so by
the above observations we now have that $H$ achieves its unique critical in
$\frak{N}_x$, as needed to prove  theorem  
\ref{tri1}.

\begin{subsection}{The Differential:  The Computation of Lemma
\ref{diff}}\label{difff}

In this section we gain our needed understanding of the
differential as expressed in lemma \ref{diff}.  
To get started note the sum in $dH$ is over all triangles but
the fact concerns only each individual one.  So we may
restrict our attention to one triangle.
One way to prove lemma \ref{diff}  
is to explicitly compute a formula for the volume
in terms of the Lobachevsky function and then find its
differential.  This method can be found carried out in
\cite{Le}.  Here we present an argument using Schlafli¹s formula for volume deformation.  This technique  has a wider
range of application as well as  being considerably more
interesting.

To start with we will recall Schlafli's formula for a
differentiable family of compact convex polyhedra  with
fixed combinatorics. Let $\{ edges \}$ denote the set of
edges and
$l(e)$ and $\theta(e)$  be the length and dihedral angle
functions associated to an edge $e$. Schlafli's formula is 
the following formula for the deformation  of the volume
with in  this family

\[ dV =- \frac{1}{2} \sum_{edges} l(e) d(\theta(e)). \]

In the finite volume case when there  are ideal
vertices the formula changes from measuring the length
of edges $l(e)$ to measuring the length of the 
cut off
edges $l^{\star}(e)$.  Let us now
recall  how
$l^{\star} (e)$ is computed.  First 
fix a horosphere at each ideal vertex. Then note from any
horosphere to a point and between any pair of
horospheres there is a unique (potentially degenerate)
geodesic segment perpendicular  to the horosphere(s).  
$l^{\star}(e)$ is the signed length of this geodesic
segment;  given a positive sign if the geodesic is out side
the horosphere(s) and a negative sign if not. 
Schlafli's formula is
independent of the horosphere  choices in this
construction, and  I will refer to this fact as the
horoball independence  principle. It is   worth recalling the
reasoning behind this principle, since the ideas involved 
will come into play at several points in what follows.  

{\bf The Horosphere Independence Reasoning:}  Recall from
the proof of observation 
\ref{fund}  that at an ideal vertex $v$  we have
the sum of the dihedral  angles satisfying 
$ \sum_{e \in v} \theta(e) =  (n-2) \pi$,
and in particular 
\[ \sum_{e \in v} d \theta(e) =  0. \]
Looking at figure \ref{horo1} we see by 
changing the horosphere at the ideal vertex $v$  that 
$l^{\star}(e)$  becomes  $l^{\star}(e) + c$ for
each $e \in v$  with $c$ a fixed constant.  Hence by our
observation  about the angle differentials
\[-2dV =   \sum_{edges} l^{\star}(e) d
\theta_e= \sum_{\{e \in v\}^c} l^{\star}(e) d \theta(e) + 
\sum_{e
\in v} (l^{\star}(e) + c )d \theta(e)  \]
and $dV$ is seen to be independent of the 
horosphere choices.

\qed

Now let us look at our prism. Let the
notation for the cut off edge lengths coincide with the
edge names in figure
\ref{clam}.
 Since we may choose any horospheres let us choose
those tangent to the hyperbolic plane which our prism is
symmetric across. In this case note the lengths of 
$(ab)^{\star}$, $(bc)^{\star}$  and $(ac)^{\star}$ are zero.
Recalling from the proof of observation \ref{fund} that  
\[ A^{\star} =  \frac{\pi + A -B -C}{2} \]
and viewing $V(d^t(x))$ as a function on 
\[ \{(A,B,C) \in (0,\pi)^3 : 0<A+B+C<\pi \} \] 
we see from Schlafli's formula that
\[ dV =  - a^{\star} dA^{\star}  - b^{\star} d B^{\star}
 - c^{\star} d C^{\star} .\]

\begin{figure*}
\vspace{.01in}
\hspace*{\fill}
\epsfysize = 2.5in  
\epsfbox{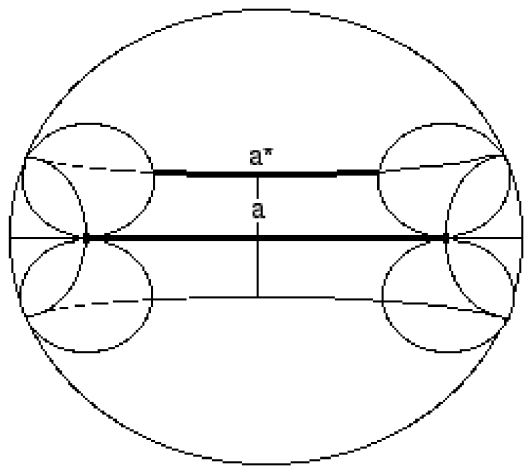}
\hspace*{\fill}
\vspace{.01in}
\caption{\label{facesl} The face of our prism containing $a$
along with the horocircle slices of the horospheres tangent
to the hyperbolic plane through which our prism
is symmetric. }
\end{figure*}

Note that lemma \ref{diff} will follow  from the
following formula.

\begin{formula} 
\[a^{\star} = 2
\ln\left(\sinh\left(\frac{a}{2}\right)\right).
\]
\end{formula}

\proof 
To begin this computation look at the face of the
prism containing
$a$
 as in figure \ref{facesl}.
 Notice this face is
decomposed into four quadrilaterals as in figure
\ref{quad}. 
Note that just as with the above reasoning
concerning the independence of horosphere choice we have
an independence of horocircle choice and
\[\frac{a^{\star}}{2} = (t^{\star} - h^{\star}) - (h^{\star}
- s^{\star}) .\]
In fact $t^{\star} - h^{\star}$  and  $h^{\star}
- s^{\star}$ are independent of this horocircle
choice as well and it is these quantities we shall
compute.  
 
\begin{figure*}
\vspace{.01in}
\hspace*{\fill}
\epsfysize = 2.5in  
\epsfbox{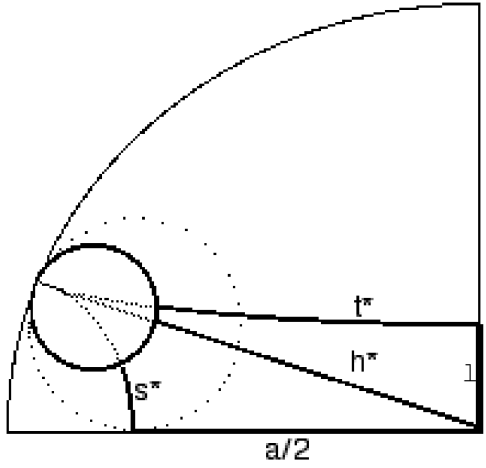}
\hspace*{\fill}
\vspace{.01in}
\caption{\label{quad}  Pictured here is one of the four 
triply right angled quadrilaterals from figure
\ref{facesl}.  Such quadrilaterals are known as Lambert
quadrilaterals and it is a well know relationship that 
$\tanh\left(\frac{a}{2}\right) = \sech(l)$, see for example
\cite{Gr}.  (In fact it follows nearly immediately from
the perhaps better known Bolyia-Lobachevsky formula.)}
\end{figure*}

Look at the figure \ref{quad} and notice using the
 horocircle tangent to the $\frac{a}{2}$ geodesic that
$h^{\star} - s^{\star} $ becomes precisely
$h^{\star}$.  Viewing this situation  as
in figure \ref{cayly} we can now read off from figure
\ref{cayly}  that
\[ h^{\star}   - s^{\star} =
-\ln\left(\sech\left(\frac{a}{2}\right)\right) .
\]

Similarly notice that 
 \[ -h^{\star}
+ t^{\star} =  \ln(\sech(l)), \]
which as observed in figure \ref{quad} implies 
\[ - h^{\star} + t^{\star}  = \ln\left(\tanh\left(\frac{a}{2}
\right)\right) .\]

\begin{figure*}
\vspace{.01in}
\hspace*{\fill}
\epsfysize = 2.5in  
\epsfbox{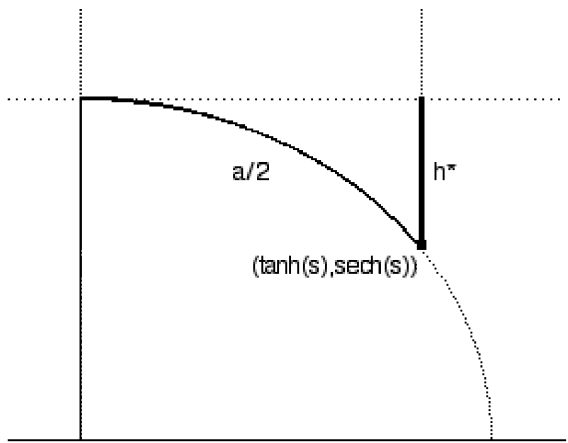}
\hspace*{\fill}
\vspace{.01in}
\caption{\label{cayly} Here we have placed the lower
triangle from figure \ref{quad}  into the upper-half
plane model sending the ideal vertex to infinity and the
$\frac{a}{2}$ segment onto the unit circle as pictured.  
Recall that the  unit  circle in this picture can be
parameterized by hyperbolic distance from
$i$ via
$\tanh(s) +  \sech(s) i$.}
\end{figure*}

With these computations we now have 
\[ a^{\star} = 2\left(
\ln\left(\tanh\left(\frac{a}{2}\right)\right) -
\ln\left(\sech\left(\frac{a}{2}\right)\right)\right)  =  2
\ln\left(\sinh\left(\frac{a}{2}\right)\right)\] 
as needed.

\qed

\end{subsection}
\begin{subsection}{Convexity:  The proof of Lemma
\ref{cony}}\label{con} 
To prove $H$ is strictly concave we start with the 
observation that  the objective function $H$  will certainly 
be a  strictly concave function on
$\frak{N}_x$ if the prism volume function $V(d^t(x))$ 
viewed as a function on   
\[ \{(A,B,C) \in (0,\pi)^3 : A + B + C < \Pi \} \]  
turned out to be strictly concave. 
In fact
it is worth noting that this implies  
$H$ is then strictly concave on all  of $(0,\pi)^{3F}$
(see section \ref{quest}). 

There are several nice methods to explore the concavity of
$V(A,B,C)$. One could simply
check  directly that
$V$'s Hessian is negative  definite (as done in 
\cite{Le}), or one could exploit the visible injectivity of
the gradient, or one could bootstrap from the concavity 
of the ideal tetrahedran's volume. It is this last method
that will be presented here.
The  crucial observation is that any
family of ideal prism can be decomposed into three ideal
tetrahedra as in figure
\ref{tet}.     So we have 
\[ V(A,B,C) = \sum_{i=1}^3 T_i(A,B,C),\]
were $T_i$ is the volume of the $i^{th}$ tetrahedra in this
decomposition.
\begin{figure*}
\vspace{.01in}
\hspace*{\fill}
\epsfysize = 2.8in 
\epsfbox{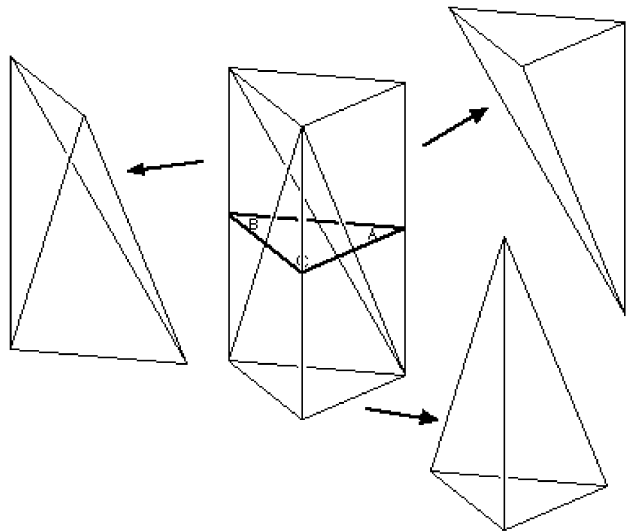}
\hspace*{\fill}
\vspace{.01in}
\caption{\label{tet} A decomposition of the ideal prism into
three ideal tetrahedra.   Notice that the
angles in this decomposition are are determined by the
affine conditions coming from the ideal vertices (see figure
\ref{horo1}) along with the condition that the angles meeting
along an edge slicing a prism face sum to $\pi$.  In particular all the
angles depend affinely on the angles
$\{A,B,C\}$.  }
\end{figure*}

Let us note some properties of the ideal tetrahedra
and its volume.  First recall from figure \ref{horo1} that the dihedral angles
corresponding to the edges meeting at a vertex of an ideal
tetrahedron are the angles of a Euclidean triangle.   In
particular the constraints 
\[ \sum_{e \in v} \theta^e = \pi \]
at each vertex guarantee that 
an ideal tetrahedron is uniquely determined by any pair of
dihedral angles $\alpha$ and
$\beta$ corresponding to a pair of edges sharing a vertex. 
Further more any pair of angles in
\[\{ (\alpha,\beta) : \alpha + \beta < \pi \}\] 
determines an ideal tetrahedron.
Note  the following fact (see \cite{Ri2}).

\begin{fact}
The ideal tetrahedron¹s volume function, $T(\alpha,\beta)$, is
strictly concave on the set 
\[\{ (\alpha,\beta) : \alpha + \beta < \pi \}\] 
and  continuous on this set's closure.  
\end{fact}

From figure \ref{tet} each of the
$\alpha_i$ and
$\beta_i$ of the $i^{th}$ tetrahedron depend on the $(A,B,C)$
affinely. So this fact immediately provides us with the
continuity assertion in lemma \ref{cony}. To exploit the
tetrahedran's concavity  we will use the following lemma.

\begin{lemma}
Let $T$ be a strictly concave function
on the convex set $U \subset \Bbb{R}^{m}$ and for each
$i$ let  $L_i$ be an affine  mapping  from $\Bbb{R}^n$ to
$\Bbb{R}^{m}$ taking the convex set $V$ into $U$. Then
the function  
\[ V(\vec{x}) = \sum_{i = 1}^{k} T(L_i(\vec{x})) \]
is strictly concave on $V$  provided $L_1 \times \dots
\times L_{k}$ is injective.
\end{lemma}

\proof
Let $l(t)$ be a line such that $l(0) = a \in V$ and 
$l(1) = b \in V$; and let $t \in (0,1)$. 
Note by the concavity of $T$ that
\[ V(l(t)) = \sum_{i=1}^3 T(L_i(l(t))) \geq  \]
 \[\sum_{i=1}^k (T(L_i(l(a)))+ t(T(L_i(l(b))) -
T(L_i(l(a)))) =
V(l(a)) + t(V(l(b)) - V(l(a))). \]
Since $T$  is in fact strictly concave the inequality
would be strict if for some $i$ we knew $L_i(l(t))$  was a
non-trivial line.  Fortunately this is guaranteed by the
injectivity of
 $L_1 \times \dots \times L_{k}$, and we are done.

\qed

Letting 
\[ L_i(A,B,C) = (\alpha_i(A,B,C),\beta_i(A,B,C))\]
we see that $V$ we will satisfy the lemma if, for example, 
 the mapping 
\[(\alpha_1(A,B,C),\alpha_2(A,B,C)),\alpha_3(A,B,C))\]
is injective. Looking
at the decomposition in figure \ref{tet} we see 
that we may in
fact choose
$\alpha_1(A,B,C) = A$,
$\alpha_2(A,B,C) = B$, and $\alpha_3(A,B,C) = C$. 
So indeed,
we have our required injectivity and  $V(A,B,C)$ is strictly
concave as needed.

\end{subsection}
\begin{subsection}{Boundary Control:  Proof of Lemma 
\ref{pushin}}\label{bound}
Before proving lemma \ref{pushin} we will rephrase it
slightly.
Namely note that the  compactness of $\overline{U}$
guarantees that 
$l(s)$ eventually hits the boundary again at some $y_1$ for some unique  $s >
0$.   So we may change the speed of our line and assume  
we are using the line connecting 
the two  boundary points, namely
\[l(s) =
    (1-s)y_0 +sy_1.\] 
So lemma \ref{pushin}  is equivalent to the following lemma.

\begin{lemma}\label{push}
For every pair of points $y_0$ and $y_1$ in $\partial \frak{N}$  
but not in $B$
with $l(s) \bigcap \frak{N} \neq \phi$ we have  
\[ \lim_{s \rightarrow 0^+} \frac{d}{ds}H(l(s))   
\mbox{  } =   \mbox{  }   \infty. \]
\end{lemma} 
{\bf Proof:}
Recalling that $H(d^t(x))  = \sum_{t \in \frak{T}} V(d^t(x))$ 
we see the lemma will follow if we demonstrate that
for any triangle 
 \[ - \infty < \lim_{s \rightarrow 0^+} \frac{d}{ds}V(d^t((s))) \leq \infty ,
\]
 and for some triangle 
\[ \lim_{s \rightarrow 0^+} \frac{d}{ds}V(d^t((s))) = \infty . \]

The boundary is expressed in terms of angle data,
so it would be nice to express the
$-2 \ln\left(\sinh\left(\frac{a}{2}\right)\right)$ coefficient in front
of the $dA^{\star}$ in
$dV$ (as computed in section \ref{difff})
in terms of the angle data.  In fact we can do even better and 
put this term in a
form conveniently decoupling  the angle and curvature.

\begin{formula} 
$ - 2 \ln\left(\sinh\left(\frac{a}{2}\right)\right) $ is equal
to 
\[
\ln(\sin(B)) + \ln(\sin(C)) - \ln \left(\frac{ \cos(A - k^t(x)) - \cos(A)
}{k^t(x)}\right)  - \ln(-k^t(x)) . \]
\end{formula}

This formula relies only on the hyperbolic law of cosines
which tells us 
\[ \cosh(a) = \frac{\cos(B) \cos(C) - \cos(A)}{\sin(A) \sin(B)}. \]
Using this relationship and the definition of $k^t(x)$  we now have
 \[ - 2 \ln\left(\sinh\left(\frac{a}{2}\right)\right)  = 
- \ln\left(\sinh^2\left(\frac{a}{2}\right)\right)   =-
\ln\left(\frac{\cosh(a) - 1}{2} \right) 
 \]   \[  = - \ln\left(\frac{\cos(B +C ) + \cos(A)}{\sin(B) \sin(C)} \right)
 = - \ln\left(\frac{- \cos(A - k^t(x)) + \cos(A)}{\sin(B) \sin(C)} \right),
\]  as needed.

Using this formula we will now enumerate the possible $y_0$ and the
behavior of $\frac{d}{ds} V(d^t((s)))$ in these various cases.  Let $C$
denote a finite constant.  We will be using the fact that if $L(s)$ is an
affine function of 
$s$ satisfying
$\lim_{s
\rightarrow 0^+}L(s) = 0$  then $\lim_{s
\rightarrow 0^+} \ln|\sin(L(l(s))|$ and the $\lim_{s
\rightarrow 0^+} \ln|D(L(s))|$  can both be expressed as 
$\lim_{s \rightarrow 0^+} \ln(s) +C$.  Furthermore for convenience let 
$d^t(y_i)   = \{A_i,B_i,C_i\}$. 

\begin{enumerate}
\item
When   $d^t(y_0)$ contains no zeros and  $k^t(y_0) \neq
0$  we have that $\lim_{s \rightarrow 0^+} \frac{d}{ds} V(d^t((s))) $ is
finite.

\item
When $d^t(y_0) = \{0,0,\pi\}$, we have 
 \[ \lim_{s \rightarrow 0^+} \frac{d}{ds} V(d^t((s))) = 
\frac{1}{2} \lim_{s \rightarrow 0^+}\ln(s)(\sigma^t(y_1 - y_0))- 
\frac{1}{2}\lim_{s
\rightarrow 0^+}\ln(s) (\sigma^t(y_1 - y_0)) +C  = C.\]

\item
In the case where   $d^t(y_0)$ contains zeros  but  $k^t(y_0) \neq
0$ for each zero (assumed to be  $A_0$ below) 
 we produce a term
in the form
 \[ \lim_{s \rightarrow 0^+} \frac{d}{ds} V(d^t((s))) = 
\lim_{s \rightarrow 0^+} \left( -\ln(s)(A_1 - A_0) \right),\]
plus some finite quantity.

\item
When $k^t(y_0) = 0$  and no angle is zero
 \[ \lim_{s \rightarrow 0^+} \frac{d}{ds} V(d^t((s))) = 
\lim_{s \rightarrow 0^+} \frac{1}{2}\ln(s)
(\sigma^t(y_1 - y_0)) +C .\]

\item
When $k^t(y_0) = 0$  and one angle, say $A_0$, in $d^t(y_0)$ is zero we have 
 \[ \lim_{s \rightarrow 0^+} \frac{d}{ds} V(d^t((s))) = 
-2\lim_{s \rightarrow 0^+}\ln(s)(A_1 - A_0)) + \lim_{s \rightarrow 0^+}
\ln(s) (\sigma^t(y_1 - y_0)) +C .\]

\end{enumerate}

So the first two cases produce finite limits.  In order to  
understand the next three  limits we make some simple observations.  First 
if  $A_0 =0$ and $l(s) \bigcap \frak{N} \neq \phi$ then $A_1-A_0 > 0$.  So
limits from the third case evaluate to $+ \infty$.     Secondly note
that   when 
$k^t(y_0) = 0$ and $l(s) \bigcap \frak{N} \neq \phi$  that $ \sigma^t(y_1 -
y_0) =A_1+B_1+C_1 -( A_0+B_0+C_0 ) < 0$ and hence the limits from the  fourth
case are $+\infty$ as well.  Combining these observations we see the
fifth case always produces a
$+\infty$ limit as well.

So for each triangle the answer is indeed finite or positive infinity.  
So all we need to do is guarantee that for some triangle we achieve
$+\infty$. To do this note that in order for $y_0$ to be on the boundary of
$\frak{N}$ and not in $B$ there is some
triangle $t$ such that   
$d^t(y_0) = \{A_0,B_0,C_0\} \neq \{0,0,\pi\}$ however, either $k^t(y_0) =0$ or
some angle is zero. 
So we have at least one triangle in case 3,4, or 5 as needed.

\qed

It is worth noting that the choice of the terminology bad for the set
$\frak{B}$ is due to the fact that at such a point   all triangles would fall
into cases one or two above, and in the process we lose our
needed control over $H$.

\end{subsection}
\end{section}

\begin{section}{The Linear Argument}\label{linn}

To see the surjectivity of $\Psi$ form $\frak{D}$ to $D$
let us assume the contrary that $\Psi(\frak{D})$ 
is strictly contained in $D$ and produce a contradiction.
With this assumption we have
 a point $p$ on the boundary of $\Psi(\frak{D})$  
inside $D$.
Note $p = \Psi(y)$ for some 
$y \in \partial \frak{D}$.  Furthermore note
$(C+y) \bigcap \frak{D}$ is empty, since otherwise for some 
$w \in C$ we would have $(y + w) \in \frak{D}$ which along
 with the fact
that $\Psi$ is an open mapping when restricted to $V$  
would force  $p = \Psi(y) =\Psi(y +w)$ to be 
in the interior of $\Psi(\frak{D})$.

At this point we need to choose a particularly nice 
conformally equivalent version of 
$y$, which requires the notion of a 
stable boundary point of $\frak{D}$.  
Before defining stability, note since  $\frak{D}$ is a convex 
set with hyperplane boundary if 
$x \in \partial \frak{D}$ such that  
$(x + C) \bigcap \frak{D} = \phi$, then 
$(x + C) \bigcap \partial \frak{D} $
is its self a convex $k$ dimensional set for some $k$.    

\begin{definition} \label{stable1}
A point in $x \in \partial \frak{D}$ is stable if 
$(x + C) \bigcap \frak{D} = \phi$  and $x$ is in the interior of 
$(x +C) \bigcap \partial \frak{D} $ as a $k$ dimensional set.
Any inequality forming  $\frak{D}$ violated in order to
make $x$ a boundary point  will be 
called a violation.
\end{definition} 

The key property of a stable point is that 
a conformal change $w \in C$ 
has $x+\epsilon w \in \overline{\frak{D}}^c$  for all 
$\epsilon > 0$ 
or for any  sufficiently small  $\epsilon > 0$ we have 
$x + \epsilon w$ must still be on $\partial{\frak{D}}$ and 
experience exactly the same violations as $x$.
  The impossibility of any other
phenomena when conformally changing a stable point is at 
the heart of the 
arguments in lemma \ref{sub1} and 
lemma \ref{sub21} below.
At this point surjectivity would follow if
for a stable  $x \in \partial {\frak{D}}$ 
we knew that  $\Psi(x)$ could not be in  $D$,
producing the  needed contradiction to our $p = \Psi(x)$
choice. 
 
We will prove this by splitting  
up the possibilities into the two cases in
lemma \ref{sub1} and lemma  \ref{sub21}.

\begin{lemma}\label{sub1}
If $x \in \partial {\frak{D}}$ is stable
and $\alpha^i(x) =0$ for $\alpha^i$ in some triangle where 
$k^t(x) <0$, then $\Psi(x)$ is not in $D$.
\end{lemma}
{\bf Proof:}
Look at an  angle slot which is zero in  triangle $t_0$ satisfying  
$k^{t_0}(x) < 0$.
View this angle as living  
between the edges $e_0$ and $e_1$.  
Note that in order for $x$ to be stable that  
either $e_1$ is a boundary edge  or 
the $\epsilon w_{e_1}$ transformation 
(with its positive side in $t_0$) 
 must be protected 
by  a zero on the  $-\epsilon$   side forcing the 
condition that 
$x +\epsilon w_{e_1} \in \bar{\frak{D}}^c$ , or else for
small enough 
$\epsilon$ we would have 
$x + \epsilon w_{e_1}$ being a conformally 
equivalent point on $\partial{\frak{D}}$  with 
fewer violations. When $e_1$ is not a boundary edge
call  this neighboring triangle $t_1$ and when it is a 
boundary edge stop this process.
If we have not stopped let  $e_2$ be another 
edge bounding a zero angle slot in $t_1$ and stop if 
it is a boundary edge.
If it is not a boundary edge  
then there are two possibilities.
If
$k^{t_1} < 0$ repeat the above procedure letting $e_1$ 
play the role of
 $e_0$ and $e_2$ the role of $e_1$ and constructing an 
$e_3$ in a 
triangle $t_2$.
If  $k^{t_1}(x) = 0$ conformally change $x$ to 
\[x +  \epsilon w_{e_1} + \epsilon w_{e_2}.\]
Notice  no triangle with 
$k^t(x)=0$ can have two zeros by lemma \ref{thing}, so
for the initial zero violation to exist 
there most be a zero on the $-\epsilon$ side of 
$\epsilon w_{e_2}$.  Once again we have determined 
an $e_3$ and $t_2$.

Using this procedure  
to make our decisions  we may continue this process  
forming a set of edges 
$\{e_i\}$ with the angle between $e_i$ and $e_{i+1}$,
 $A^{i,i+1}(x)$, always 
equal to zero.   I'll call such a set an accordion, see in figure
\ref{sna}  for an example. 
\begin{figure*}
\vspace{.01in}
\hspace*{\fill}
\epsfysize = 2in  
\epsfbox{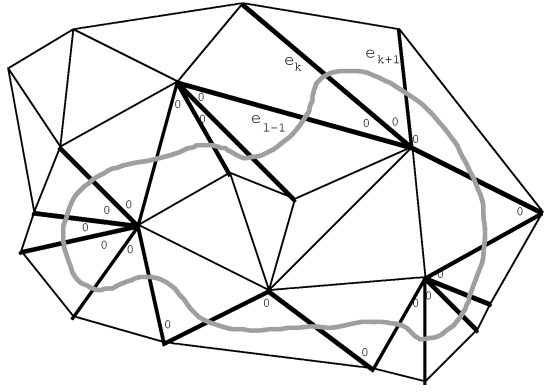}
\hspace*{\fill}
\vspace{.01in}
\caption{\label{sna}  Here we have an accordion,
the set of edges forming a loop indicated  with the squiggly
line. Notice to geometrically realize the 
 zero angles in the pictures corresponds to squeezing the 
accordion.  The algebraic inability to squeeze various accordions is at the heart
of the $\Psi$'s surjectivity.
 }
\end{figure*} 
Since there are a finite number of edges 
an accordion either stops at a boundary edge or 
the accordion is endless.  If the accordion is endless then eventually the 
sequence $\{e_i\}_1^{\infty}$  will have some $k<l$ such
that $e_k = e_l$  and $e_{k+1} =e_{l+1}$,  see figure
\ref{sna}. (This is true by the pigeon hole principle since 
some edge $e$ will appear an infinite number of 
times in this 
list and among its infinite neighbors there must
 be a repeat). 
We can produce a contradiction to this occurring. 
To do it first note if  $e_{i}$ and 
$e_{i+1}$ are in $t_i$ then
$A^{i,i+1}(x) = \psi_{t_i}^{e_i} + \psi_{t_i}^{e_{i+1}}$.  So for 
the set of edges $\{e_i\}_{i=k}^{l-1}$ we have 
\[ 0 = \sum_{i = k}^{l-1} A^{i,i+1}  = 
\sum_{i=k}^{l-1} \psi^{e_i}(x) > 0 \]
our needed contradiction. 

In the case the sequence did hit the boundary 
perform the accordion construction  in the opposite direction.  
If we don't 
stop in this direction we arrive at the same 
contradiction.  If we did then
this computation still produces a contradiction 
on the accordion  with the 
two boundary edges, since for a boundary edge 
in the triangle $t$ we 
have
 $\psi^e_t(x) = \psi^e(x) \in \left(0,\frac{\pi}{2}\right)$.

\qed

\begin{lemma}\label{sub21}
If a stable $x$ satisfies the condition that 
if $\alpha^i(x) =0$ then $\alpha_i$ is in a 
triangle $t$ 
with
$k^t(x) = 0$, 
then $\Psi(x)$ is not in $D$. 
\end{lemma}
{\bf Proof:}
In this case, in order for $x$ to be a boundary 
point of $\frak{D}$ 
for some $t$ we have that 
$k^t=0$. We will be looking at the nonempty set of 
all triangles with $k^t = 0$, $Z$.  
The first observation needed about $Z$ is that it is
not all of $M$ and has a non-empty internal boundary 
(meaning $\partial Z - \partial M$).
To see this note 
\[\sum_{t \in \frak{T}} k^t(x) = \sum_{e \in v} A^i -\pi F 
 = \pi {\partial} V + 2 \pi (V - \partial V) -  \pi F \]
\[  =
2 \pi V - (\pi {\partial}  V+ 3 \pi F) +  2 \pi F  
=   2 \pi V - 2 \pi E + 2 \pi F  = 2 \pi \chi(M) < 0,\]
so there is negative curvature somewhere.

By the stability  
of $x$ once again 
there can be no conformal transformation capable of 
moving negative 
curvature into this set. Suppose we are at an internal
 boundary $e_0$
edge of $Z$, call the triangle  on the $Z$ side of 
the boundary edge $t_0$ and the triangle on the
non-boundary edge  $t_{-1}$.      Since  $t_{-1}$ has 
negative 
curvature the obstruction to 
the $\epsilon w_{e_0}$ transformation
being able to move 
curvature out of $Z$ must be due to  $t_0$.  
In order for $t_0$ to protect against 
this there must be zero along $e_0$ on the $t_0$ side.

Now we will continue the attempt to suck
curvature out with a curvature vacuum.
Such a vacuum is an element of $C$ indexed by a set of 
$Z$ edges.  
The key observation in forming this vacuum is once again
lemma 
\ref{thing} telling us
if an angle in $t$ is zero and $k^t(x) = 0$ 
then there is only one zero angle in $t$.  
Let $e_1$ be the other edge sharing the unique zero angle
 along 
$e_0$ in $t_0$ and if 
$e_1$ is another boundary edge we stop. If $e_1$ is not a
 boundary edge
 use $\epsilon(w_{e_1} + w_{e_0})$ 
to continue the effort to remove curvature.  Continuing this 
process forms a completely determined
set of edges  and triangles, $\{e_i\}$ and $\{t_i\}$,   and 
a sequence of conformal   
transformations $\epsilon \sum_{i =0}^n w_{e_i} \in C$, see figure \ref{vac}.
\begin{figure*}
\vspace{.01in}
\hspace*{\fill}
\epsfysize = 2in  
\epsfbox{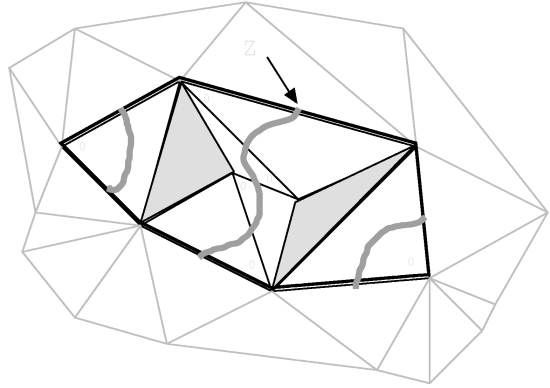}
\hspace*{\fill}
\vspace{.01in}
\caption{\label{vac} In this figure the region of zero curvature, $Z$,  
is the set with the bold boundary as labeled.   We have also attached to $Z$
all its vacuums and indicated them with the squiggly lines.  As we can see in 
the figure these vacuum transformations fail  to suck curvature out of
$Z$, just  as in the proof of lemma \ref{sub21}.  Notice the set  $S$,
$Z$ minus its vacuums,  is in this case the pair of shaded
triangles. }
\end{figure*} 

We will now get some control over this vacuum.  Note a
vacuum never 
hits itself since if there is a first pair 
$k <l$ such that
$t_k =t_l$ then $t_k$ would have to have two zeros 
and zero curvature, which lemma \ref{thing} assures us is 
impossible.  So any vacuum hits a boundary edge or  
pokes  through $Z$ into $Z^c$.

In fact with this argument we can arrive at 
the considerably stronger fact 
that two vacuums can never 
even share an edge.  
To see this, call a  vacuum's side
boundary any edge  of a triangle in the vacuum 
facing a zero.  
Now simply note if the intersection of  two vacuums 
contains an edge then
it contains a first edge $e_i$ with respect to one of 
the vacuums.
There are two possibilities for this edge. One is that 
$t_{i+1}$ has two zeros 
and $k^t(x) =0$, which we showed was impossible in
 the previous paragraph.
The other is that $e_i$ is a side boundary 
of both  vacuums.
In this case we have an edge facing zero angles 
in both directions in triangles with zero curvature,  so this
would force $\psi^e(x) = \pi$, a contradiction.  So either
case  is impossible, and indeed no distinct vacuums share an
edge. 
 
Let 
$S$ be the removal
from $Z$ of all these vacuums, see figure \ref{vac}.  
First I'd like to note that
$S$ is non-empty.  Note every vacuum has  
side boundary.  Since vacuums cannot intersect themselves or 
share edges with distinct 
vacuums,  $S$ 
would be nonempty if side boundary had to be 
in $Z$'s interior.  
Look at any side boundary edge $e$ of a fixed vacuum. 
 Note $e$ 
cannot be on $\partial Z - \partial M$ since then the vacuum 
triangle it belonged to 
would have at least two zeros  and $k^t(x) =0$.  
Furthermore, 
$e$ cannot be on $\partial M$ 
since then $\psi^e(x) = \frac{\pi}{2}$.  
So indeed $S$ is nonempty.

Now let's observe the following formula.
\begin{formula}\label{set}
Given a set of triangles $S$ 
\[ \sum_{\{e \in S \}} \theta^e(x) =
\sum_{t \in S} \left( \pi  - \frac{k^t(x)}{2}\right) 
+ \sum_{e \in \partial S-\partial M}  
\left(\frac{\pi}{2} - \psi^e_{t}(x)\right) , \]
with the   $t$ in the $\psi^e_{t}(x)$ term being the triangle on 
the  non-$S$ side of $e$.   
\end{formula}
{\bf Proof:}
\[ \sum_{\{e \in S \}} \theta^e(x) = 
\sum_{e \in S} (\pi - \psi^e(x)) \]
\[=
\sum_{e \in S-\partial M} \left( \left(\frac{\pi}{2} - \psi_{t_1}^e(x)\right) 
+ \left(\frac{\pi}{2} - \psi_{t_2}^e(x)\right)\right)  
+\sum_{e \in \partial M \bigcap S} 
\left(\frac{\pi}{2} - \psi_{t}^e(x)\right) \] 
\[ = \sum_{t \in S} \left( \pi  + \frac{ \pi - l^t(x)}{2}\right) 
+ \sum_{e \in \partial S -\partial M}  
\left(\frac{\pi}{2} - \psi^e_{t}(x)\right)  \]
\[ =\sum_{t \in S} \left( \pi  - \frac{k^t(x)}{2}\right) 
+ \sum_{e \in \partial S-\partial M}  
\left(\frac{\pi}{2} - \psi^e_{t}(x)\right). \]

 \qed

Now every edge in $\partial S -\partial M$ faces a zero on its $S^c$ side
in a triangle with $k^t(x) =0$ (see figure \ref{vac} once again), so
\[ \sum_{e \in \partial S-\partial M}  
\left(\frac{\pi}{2} - \psi^e_{t}(x)\right) =0.\]
Similarly each triangle has zero curvature so from the above 
formula we have 
\[   \sum_{\{e \in S \}} \theta^e(x) = |S| \pi \]
violating condition $(n_2)$.  So we have constructed a violation to 
$(n_2)$ and $\Psi(x)$ cannot be in $D$ as need.

\qed

\end{section}

\begin{section}{Generalizations and Comments}\label{itat}
\begin{subsection}{Corners, Cones and Cusps}

Its worth noting that the proof  of theorem \ref{tri1}
relies in no way on the restriction that $p_v(x) = 2 \pi $ at internal vertices
and $p^v(x) = \pi$ at boundary vertices.  This assumption  gives the simplest 
case, of producing hyperbolic surfaces with geodesic boundary.   If this
condition is relaxed to using any positive numbers, the same exact proof goes
through to produce triangulations of surfaces with cornered boundary  and cone
singularities.  

\begin{figure*}
\vspace{.01in}
\hspace*{\fill}
\epsfysize = 2in 
\epsfbox{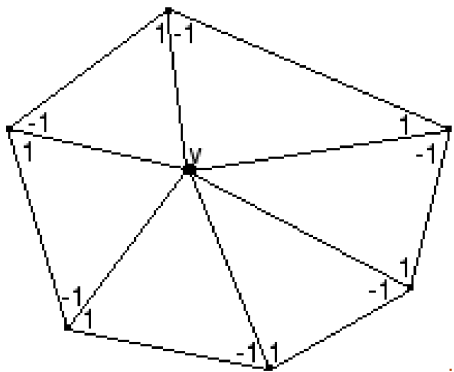}
\hspace*{\fill}
\vspace{.01in}
\caption{\label{flowt} Let $v \in V_{\infty}$ and call  the vector in 
this figure $f_v$.  Notice that
this vector is $\sum_{e \in v}  w_e$ (with the
correct sign choices) hence a conformal transformation.
For the purpose here, the most   important  observation is that  
constant multiples of
this vector are precisely the  conformal transformations involving $v$ 
that preserve the condition that the angles at $v$ remain zero.   In general
$span_{v \in V} \{f_v\}$ arises  naturally as the  set of   
transformations preserving each triangle¹s curvature.}
\end{figure*}

By far the most interesting case is when this condition  is set to
$p_v(x) = 0 $ at some vertices,  and we produce surfaces with cusps.  For the
discussion here  let us suppose that this
condition is placed at a set of interior vertices denoted $V_{\infty}$ 
and that such vertices are
isolated in the sense that no single triangle has two of its vertices in
$V_{\infty}$.  Notice $p_v(x) =0$ forces us to make all the angles at $v$
identically zero.   Letting $E_{\infty}$ and $F_{\infty}$  
 denote the edges and faces with a vertex in $V_{\infty}$ respectively,
our space of angle systems, $\frak{N}$   
is restricted to $(0,\pi)^{3F - E_{\infty}}$.  The  conformal
transformations are spanned  
by the $w_e$ at edges in  $E - E_{\infty}$
along with the vectors $f_v$ from figure \ref{flowt} for each $v \in
V_{\infty}$.

On our new $\frak{N}$  we can put on the
same objective function $H$, though for faces in $F_{\infty}$ the prism
degenerates to the polyhedron in figure \ref{deg}.
\begin{figure*}
\vspace{.01in}
\hspace*{\fill}
\epsfysize = 2in 
\epsfbox{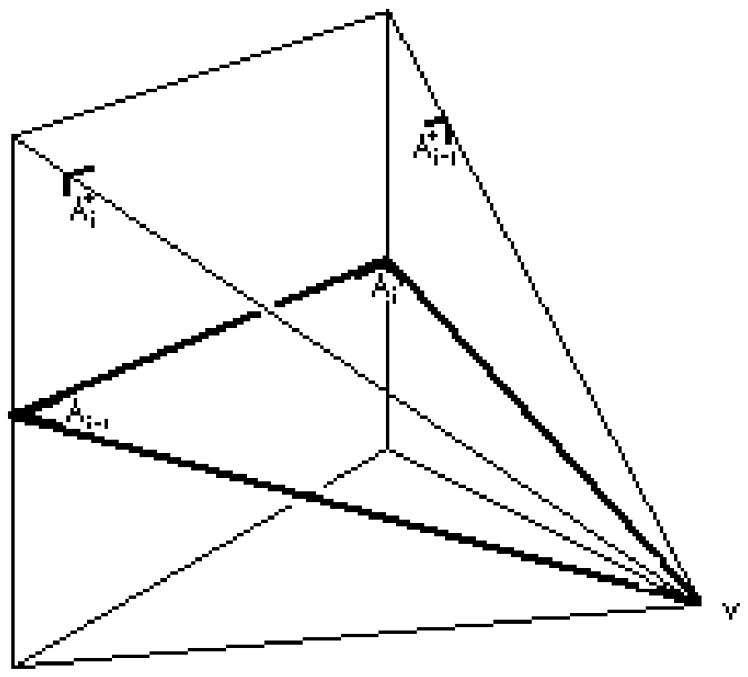}
\hspace*{\fill}
\vspace{.01in}
\caption{\label{deg} Here is the notation for our degenerate prism, when $d =
\{A_{i-1},A_{i},0\}$ with
$A_k$ greater than zero.  Notice the  side of our original prism where the
zero angle lives degenerates to zero length, and the pictured vertex $v$ is on
$S^{\infty}$.}
\end{figure*}
The same exact same arguments as in section \ref{prot} allow us to see that at a
critical point the edges in $E -E_{\infty}$ fit together.  The remaining edges
are all infinite and so we can certainly glue the triangles of a flower at $v \in
V_{\infty}$ together.   To understand this situation, pick a cyclic order on
the triangles this  flower 
$\{t_1, \dots , t_n\}$ and denote the ordered pair of $E_{\infty}$  edges  of
$t_i$  as
$\{e_{i-1}, e_{i}\}$.  With this notation place the realizations of the
$d^{t_i}(x)$  in the upper half plane as in figure
\ref{fit}.

\begin{figure*}
\vspace{.01in}
\hspace*{\fill}
\epsfysize = 2in 
\epsfbox{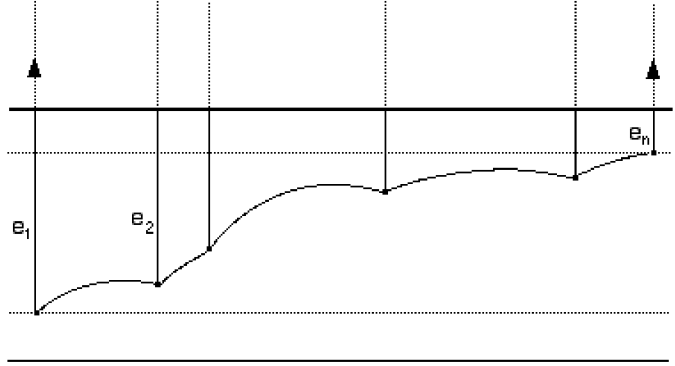}
\hspace*{\fill}
\vspace{.01in}
\caption{\label{fit}
In this figure we have placed the realizations of the triangles in our flower
next to each other.   When we glue the realization of   $e_0$ to
the realization of  $e_{n}$ we  will typically  not get a cusp, in fact 
such a gluing will usually produce  a non-complete hyperbolic surface (see \cite{Th}). 
However, when the dotted horocircles in the figure  agree,  we may use an
isometry represented by a Euclidean translation of the upper-half plane to
identify $e_0$ and $e_n$ and produce a cusp with piecewise geodesic
boundary.  The thick line is a
horocircle which could be used to simultaneously cut all the infinite
triangle sides.}
\end{figure*}

At this point we need to know that at a critical point we satisfy the cusp
producing condition as in figure \ref{fit}; the trick to accomplishing this is to use
the horocircle independence principle introduced in section \ref{poly} along with
a "holonomy" argument.    Notice that
$l^{e_i}(u)$ is infinite.  As usual by using a horoball cut off the
length of 
$(l^{e_i}(u) - l^{e_{i-1}}(u))$ is well defined.  
 This observation being  true for each
$t_i$ independently allows us to use a simultaneous cut off of the whole
realization in the upper-half plane as in figure
\ref{fit}. 

Using the notation of figure \ref{clev}  we see that $e_0$ and $e_n$ matching up
correctly is equivalent to 
\[    e_{n+1}^{\star} - e_0^{\star} = \sum_{i =1}^n
(e_i^{\star} - e_{i-1}^{\star}) =0 .\]
(This observation is what I have referred to as an holonomy argument and it
originally came up in this type of  situation in  Bragger's treatment of the
Euclidean case, see
\cite{Be}). 

The usual angle formulas (from the proof of observation \ref{fund})  hold for
$A_i^{\star}$ and
$A_{i-1}^{\star}$,  for exactly the same reasons.  
Just as in section \ref{difff} we will be computing $dV$ via Schlafli's formula 
and will use any horosphere to cut off the vertex corresponding to our $v \in
V_{\infty}$ and the horospheres tangent to the specified
$H^2$ to cut off the remaining "prism" vertices. 
Using the notation of figure \ref{clev} at a
critical point we have
\[ 0 = dH (f_v) = \sum_{e \in v} (e^{\star \star}_{i} - e^{\star \star}_{i}).
\] However we can also see  in  figure \ref{clev}
that $e^{\star} = e^{\star \star} +
\ln(2)$, so this equation also implies the needed equation in the
previous paragraph. Hence theorem \ref{tri1}  holds in the cusp case as
well.

\begin{figure*}
\vspace{.01in}
\hspace*{\fill}
\epsfysize = 2in 
\epsfbox{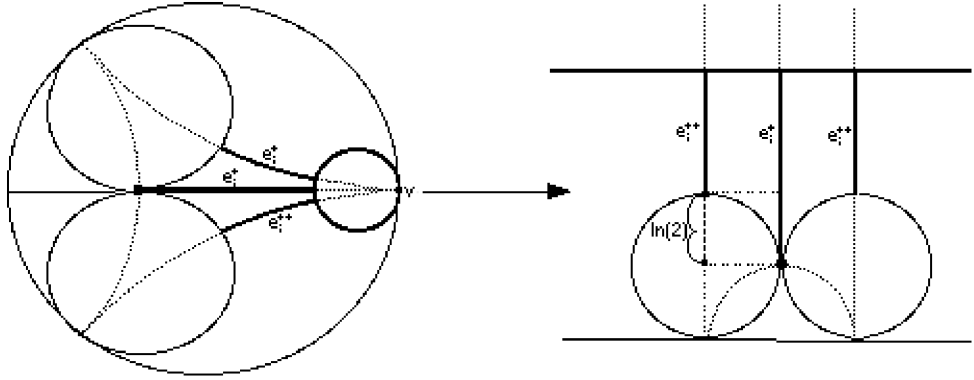}
\hspace*{\fill}
\vspace{.01in}
\caption{\label{clev}
The face  of the degenerate prism containing the edge
$e_i$  is pictured, along with its horosphere cuts.  For
convenience
$e_i^{\star}$ will denote both the cut off version of $e_i$'s realization with
respect to the critical point's data  and $e^{\star}$'s  length. Similarly for
$e_i^{\star
\star}$, the degenerate prism's labeled edge. Notice by transforming this picture
to the upper-half  plane model we can immediately compute that
$e^{\star} = e^{\star \star} + \ln(2)$.
}

\end{figure*}

\end{subsection}

\begin{subsection}{The Convex Case}
Theorem \ref{cir11} can of course be extended to the above mentioned cases as
well. Much more importantly though, it can be extended
to the entire convex case where $\theta^e(p) \in (0,\pi]$.
This generalization is very useful to produce as
corollaries other versions of the Thurston-Andreev theorem.
The price is that there are new affine conditions
placed on $p \in \Bbb{R}^E$ and the analogs to the linear
arguments in section 
\ref{linn} become considerably more complicated.  
I'll state the theorem whose careful proof can be found in \cite{Le}. In
order to articulate the new conditions we  need certain 
 snake and a loop concepts in a triangular decomposition.

\begin{definition}\label{snakes}
  A snake is a finite
directed sequence of edges $\{e_i\}_{i=k}^{l}$
directed in the following sense: if $k <l$ we 
start with the edge $e_k$ between $t_{k-1}$ and $t_k$,  
then we require $e_{k+1}$ to be one of the remaining edges on $t_k$.  
Then letting 
$t_{k+1}$ be the other face associated to $e_{k+1}$ we require $e_{k+2}$ 
to be 
one of the other edges of $t_{k+1}$ and so on until some 
tail edge $e_l$ 
and tail face $t_l$ are  reached, and 
if $l< k$ we reverse the procedure and subtract   from rather than add to
the index.  
A loop is a  snake $\{e_i\}_{i=l}^{k}$ 
where   
$e_k = e_l$ and $t_k = t_l$.
\end{definition}

It is a condition on  snakes and loops 
which allows  one to articulate
the remaining necessary conditions. It should be 
clear already how such objects
are naturally born  from  the arguments in lemmas
\ref{sub1} and \ref{sub21},in fact  one can see the accordion in figure \ref{sna}
and the vacuums in figure \ref{vac} for examples of a loop and snakes
respectively.   As defined  snakes and loops can self intersect  and
there will be an infinite number of such objects,  and it  is  worthwhile to first isolate
a  finite sub-set that does the job.

\begin{definition}\label{barbell}
A set of edges $\{e_i\}_{i=k}^l$  is called embedded if 
$e_i \neq e_j$. A snake $\{e_i\}^l_k$  is said to double back on itself 
if we have a pair of non-empty sub-snakes with
$\{e_i\}_{m}^{n}$ and $\{e_i\}_{k-m}^{k-n}$ containing the same edges.  
A barbell is a loop which doubles back on itself and such that 
$\{e_i\}_{i=k}^l / \{e_i\}_{i=m}^n $ is embedded.  A balloon  is a snake 
which doubles back on itself with 
$\{e_i\}_{i=k}^l / \{e_i\}_{i=m}^n $ embedded and such that $e_l = e_k$.
\end{definition} 

With this terminology the remaining necessary conditions are  
\begin{eqnarray*}
(n_3) & \sum_{i=k}^{l-1} \theta^e_i(p) < |k-l| \pi & 
\mbox{ when } \{e_i\}_{i=k}^l \mbox{ is an embedded loop or barbell,}
\end{eqnarray*}  
and 
\begin{eqnarray*}
(n_4) & \sum_{i =k}^{l} \theta^e_i(p) < (|k-l| +1)  \pi & 
\mbox{ when } \{e_i\}_{i = k}^l \mbox{ is an embedded sake or balloon.}
\end{eqnarray*}

With these conditions we can completely characterize the angles arising in convex
ideal disk patterns.

\begin{theorem}\label{cir1}
That $p \in \Bbb (0,\pi]^E$ and satisfies $n_i$ for each $i$ 
is  necessary and sufficient 
for $p$ to be equal to $\Psi(u)$ for some uniform angle system $u$.
Furthermore this $u$ is unique. 
\end{theorem}

\end{subsection}

\begin{subsection}{Some Natural Questions}\label{quest}
Here I list four cases where I think it would be particularly nice to attempt
to apply these  hyperbolic volume techniques.

\begin{enumerate}

\item {\bf The Spherical Case}
The use of hyperbolic volume in this paper to
solve theorems
\ref{tri1} and
\ref{cir11} could conceivably be used to prove the analogous questions in the
spherical case. The polyhedron to be used now becomes the twisted prism in figure
\ref{sphere} (also observed as the right object for this game by Peter Doyle). 
It is easy to see that the
critical points of the hyperbolic volume objective function are once again
precisely the  uniform angle system; but the objective function fails to be
convex. Can this objective function still be used to arrive at the analogous
results? (See
\cite{Le} for a more detailed account of this situation.)

\begin{figure*}
\vspace{.01in}
\hspace*{\fill}
\epsfysize = 2in 
\epsfbox{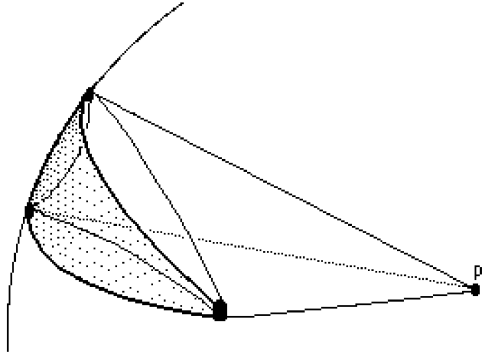}
\hspace*{\fill}
\vspace{.01in}
\caption{\label{sphere}
Fix a point $p$ in $H^3$ viewed in the Poincare Model.  Imagine
 $S^{\infty}$ is given the unit sphere's induced metric from $\Bbb{R}^3$.  
Given a geodesic triangle on with a hyperbolic
geodesic to $p$. The polyhedron of interest  is the  convex hull of these
geodesic segments.  }
\end{figure*}

\item{\bf The Non-compact Case}

Analogs to theorems \ref{cir11} and \ref{tri1} also exist  in non-compact situations.  
Let $\frak{T}$ be a locally finite triangulation on a topological surface
and let  $S$ be some set of  triangles  in $\frak{T}$.   Given any finite  $\hat{S} \subset S$ 
let $\partial_{S}\hat{S}$  be the subset set of $\hat{S}$  boundary edges 
which are not $S$ boundary edges. 
Denote as  $(\hat{n}_2)$ the condition that for any  set $S$ there is some 
 $\hat{S} \subset S$ such that 
 \[ \sum_{e \in \hat{S}}  \theta^e(p)  - \pi |\hat{S} |  >  \pi |\partial_S \hat{S}|. \]
Note  this condition is  equivalent to $(n_2)$ when $\frak{T}$ is finite.

 Let condition  $(\hat{n}_4)$  be that for any {\bf infinite}  snake $\{e_i\}$ 
there exist $N_i$   such that
\[ \sum_{i  = N_1}^{N_2} \psi^{e_i}(p)  >  \pi   .\]
Note this  condition is  automatically satisfied in the finite case as well.

Let a  geodesic  triangulation  satisfy having all its  angles in $(\epsilon,\pi]$ 
and all its $k^t $ in $ [-\pi, -\epsilon) $ for some $\epsilon >0$.
Then it  is straightforward  to verify that $(\hat{n}_2)$ and $(\hat{n}_4)$ are in fact necessary.
  The following corollary to theorem \ref{cir11}  guarantees these condition are always
sufficient.

\begin{corollary}
Given  $p  \in  (0,\pi)^{(E - \partial E)} \times \left(0,\frac{\pi}{2}\right)^{\partial E}$  
satisfying  
$(n_1)$ ,  $(\hat{n}_2)$,    and  $(\hat{n}_4)$,  
then $p$ is  an ideal disk pattern.
\end{corollary}

Two question immediately arise.  First  would be nice to know conditions guaranteeing when
there is a solution forming a complete surface.   Secondly it would be nice to to understand the
uniqueness, or more likely the moduli of space of possible solutions. 

\item{\bf The Topological Case}
Given a topological triangular decomposition $\frak{T}$ theorem \ref{cir11} tells
whether a given set of intersection angle data can be geometrically realized as an Delaunay
ideal disk pattern. It would be nice to answer: for which $\frak{T}$
does consistent  intersection angle data exist.  There are some several easy answers.  For example   letting $\{n_1,n_2,n_3\}$ be  the degrees of the vertices of any triangle, if  we have
\[ \frac{1}{n_1} +  \frac{1}{n_2} + \frac{1}{n_2} < \frac{1}{2} \]
and 
\[ 0  <  \frac{1}{n_1} +  \frac{1}{n_2} - \frac{1}{n_2}, \]
then the angles $\frac{2 \pi}{n_i}$ satisfies  theorem \ref{tri1}  (see figure \ref{xamp}).  
It would be nice to find an identification procedure that worked for any triangulation.


Furthermore as observed in section \ref{con}  $H$ is convex through
out $\frak{N}$. A critical point here is the maximal volume associated
polyhedron among all possible realizations.  Such a critical point turns out to
be  a particularly symmetric triangulation, where  triangles $t_1$ and $t_2$ 
sharing an edge $e$ with edges given by 
$\{a_i,b_i,e\}$ will satisfies  $a_1^{\star} + b_1^{\star} = a_2^{\star} +
b_2^{\star}$.  In the presence of such a critical point we then have a canonical realization of $\frak{T}$,
and it would be nice know for which $\frak{T}$ does this realization exists?

\begin{figure*}
\vspace{.01in}
\hspace*{\fill}
\epsfysize = 2.5in 
\epsfbox{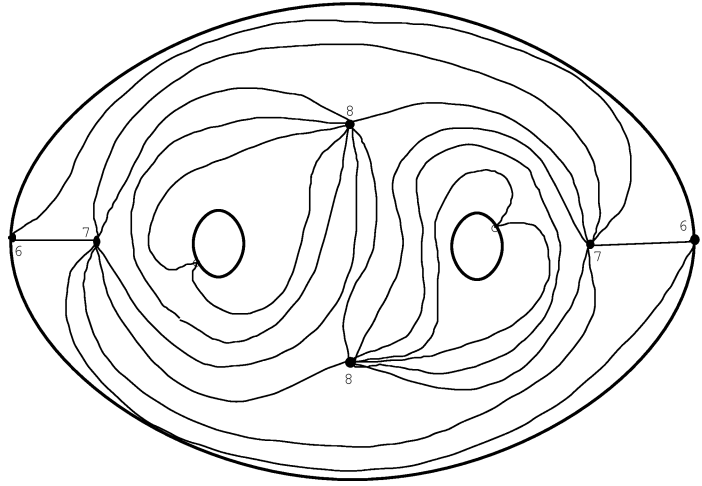}
\hspace*{\fill}
\vspace{.01in}
\caption{\label{xamp}
A topological  triangulation guaranteed to have a geometric realization as
a Delaunay triangulation.  This guarantee still hold when we glue together copies of the pictured pairs of pants.  }
\end{figure*}
 
\item
{\bf The Non-ideal Case} It would be nice to extend this use of hyperbolic volume
to the non-ideal cases, and in particular prove theorems \ref{tri1} and \ref{cir11}
in these cases.

In the sub-ideal case (when the
vertices of the hyperbolic polyhedra associated to the ideal disk pattern are finite) 
the natural hyperbolic objective function created using the corresponding  finite prisms can once again by Schlafli¹s formula be seen to have uniform critical points in a "conformal class".
However to  get started one must first construct the perpendicular edge lengths  corresponding to  each vertex (the $(ab)^{\star}$,$(ac)^{\star}$, and  $(bc)^{\star}$ in figure \ref{clam}),    so that the angles sum to $2 \pi$ at  each vertex.  In the process the objective function becomes de-localized, making boundary control and concavity difficult to verify,  and the corresponding the theorem \ref{tri1} difficult to get a handle on. 

It is worth noting one can use theorem \ref{cir1}  to achieve
 disk packing  and super-ideal versions of theorems \ref{cir1}, but one needs to
retriangulate a bit and the process feels a bit synthetic.  It would be nice to have direct hyperbolic volume techniques in this case as well.

\end{enumerate}
\end{subsection}

\end{section}

\bibliographystyle{plain}
\bibliography{thur}

\end{document}